\documentclass{article}
\usepackage{graphicx}
 \usepackage{mathptmx}
\usepackage{amsmath, amstext,amssymb,amsfonts}
\usepackage[english]{babel}
\usepackage{delarray}
\usepackage{mathptmx}
\usepackage{amsmath, amstext,amssymb,amsfonts}
\usepackage[english]{babel}
\usepackage{delarray}

\DeclareFontFamily{U}{mathx}{\hyphenchar\font45}
\DeclareFontShape{U}{mathx}{m}{n}{
      <5> <6> <7> <8> <9> <10>
      <10.95> <12> <14.4> <17.28> <20.74> <24.88>
      mathx10
      }{}
\DeclareSymbolFont{mathx}{U}{mathx}{m}{n}
\DeclareFontSubstitution{U}{mathx}{m}{n}
\DeclareMathAccent{\widecheck}{0}{mathx}{"71}
\DeclareMathAccent{\wideparen}{0}{mathx}{"75}

\numberwithin{equation}{section}
\newcommand{\R}{\mathbb R}    \newcommand{\Co}{\mathbb C} \newcommand{\Cop}{{\mathbb C}^{+}}   \newcommand{\Cod}{{\mathbb C}^{-}}
                          
 \newtheorem{thm}{Theorem} \newtheorem{lem}{Lemma}   

\date{}
\title{A note on analytic continuation   of characteristic functions}
\author{\large{ Saulius  Norvidas}}
\date{\footnotesize Institute of Data Science and Digital Technologies, Vilnius University, \\ Akademijos str. 4, Vilnius LT-04812, Lithuania\\
 ({\rm{e-mail: norvidas{@}gmail.com}})}
\begin{document}

\maketitle
 {{ {\bf Abstract}}}
 We derive necessary and sufficient conditions for a continuous bounded function $f: \R\to \Co$ to be a characteristic function of a probability
measure. The Cauchy transform $K_f$ of $f$ is used as  analytic continuation of $f$ to the upper and lower half-planes in $\Co$. The conditions
depend on the behavior of  $K_f(z)$ and its derivatives on the imaginary axis in $\Co$. The main results are given
 in terms of  completely monotonic and absolutely  monotonic functions.
 
{\bf Keywords}:  Characteristic function; complex-valued  harmonic function;   Cauchy transform;  analytic function; completely monotonic function;  absolutely  monotonic function. 

{\bf  Mathematics Subject Classification}:   30E20 - 60E10

\section{ Introduction }
Suppose that  $\sigma$ is a probability measure on $\mathbb{R}$. Let
\begin{equation}\label{1.1}
\widehat{\sigma}(x)=\int_{-\infty}^{\infty}e^{-ixt}\,d\sigma(t)
\end{equation}
be  the Fourier transform of $\sigma$. In the language of probability theory, $f(x):= \widehat{\sigma}(-x)$ is called the  characteristic function of  $\sigma$ (see [4, p. 10]).   Any characteristic function $f$ is continuous  on $\mathbb{R}$ and   satisfies   $f(-x)={\overline{f(x)}}$\ for all $x\in\R$.  In particular, such an $f$ is real-valued if and only if  it is the Fourier transform of a symmetric distribution $\sigma$ [4, p. 30], i.e. if $\sigma$ satisfies  $\sigma(-A)=\sigma(A)$ for any measurable $A\subset \mathbb{R}$.

A   function $u=u_1+i u_2$ in a domain $D\subset\Co$ is called complex-valued harmonic if both $u_1$  and $u_2$ are real harmonic functions in $D$. This means that  $u_1$ and $u_2$ are twice continuously differentiable on $D$ and satisfy there   $\Delta\ u_1= \Delta\ u_2=0$, where $\Delta$  is the Laplace operator
 \[
 \Delta=\frac{\partial^2}{\partial x^2}+\frac{\partial^2}{\partial y^2} .
\]

 For a  function $\varphi: \R\to\Co$, the Dirichlet problem on the complex upper half-plane $\Co^{+}=\{ z=x+iy\in\Co: \ y>0\}$ is to extend $\varphi$ on $\Co^{+}$ to a complex-valued harmonic function $u_{\varphi}$, so that $u_{\varphi}(z)$ tends to $\varphi(x_0)$ as $z\in\Cop$ tends to $x_0$ for each $x_0\in\R$. Let us remember one of the most important  case of this problem. If   $\varphi$ is a bounded continuous function on $\R$, then
\begin{equation}\label{1.2}
u_{\varphi}(x,y)=\frac1{\pi}\int^{\infty}_{-\infty}\frac{y}{(x-t)^2+y^2}\,\varphi(t)\,dt=\Bigl(P_y\ast  \varphi\Bigr)(x),
\end{equation}
 is the unique solution to the Dirichlet  problem on $\Cop$ that is bounded in $\Cop$. Here
\begin{equation}\label{1.3}
P_y(x)=\frac{1}{\pi}\frac{y}{x^2+y^2},\quad y>0,
\end{equation}
$x\in\R$,   is called the Poisson kernel for $\Cop$.

In this paper, we will consider the problem of characterizing of conditions for $\varphi$ to be a characteristic function in terms of its analytic continuation in $\Co\setminus\R$.    If such a $\varphi$ corresponds to a  symmetric distribution on $\R$, then a similar problem  for the  harmonic continuation (1.2)   was solved in \cite{No}. Here and for later use  we  need the
notion of completely monotonic function (see [7, p.p. 144-145]).
 A function $\omega: (a, b)\to \R$,  $-\infty\le a<b\le\infty$,  is said to be  completely monotonic if it is infinitely differentiable  and
\begin{equation}\label{1.4}
(-1)^n \omega^{(n)}(x)\ge 0
\end{equation}
for each $x\in(a, b)$ and all  $n=0,1,2,\dots$. A function $\omega: [a, b]\to \R$  is called completely monotonic   on $[a, b]$ if it is there continuous and completely monotonic   on $(a, b)$.

\begin{thm}\label{t1}{\rm{[5,  Theorem  2].}}
Suppose that $\varphi: \R\to\R $ is a   bounded continuous  even function   and  $\varphi(0)=1$. Then $\varphi$ is a characteristic function if and only if the function
\[
y\to u_{\varphi}(0,y)
\]
 is completely monotonic on $[0,\infty)$.
\end{thm}

Other form  of this theorem (in other terms) has been shown by Egorov \cite{Eg}. Moreover,  \cite{Eg} deals with absolutely integrable and infinitely differentiable  $\varphi$ satisfying more other conditions.

 If $\varphi$ is a complex-valued function, then it is easy to see that Theorem 1 fails   in general. Indeed, suppose that
\begin{equation}\label{1.5}
\varphi(x)=(1-\alpha)e^{-ix}+\alpha e^{ix},
\end{equation}
where $ \alpha\in\R$. Since
\[
\widehat{P_y}(\xi)= e^{-y|\xi|}, \ \xi\in\R,
\]
for all $y>0$,  we obtain   by straightforward calculation in (1.2) that  $u_{\varphi}(0,y)= e^{-y}$. Therefore,     $u_{\varphi}(0,y)$ is completely monotonic on $[0,\infty)$ for all  $\alpha\in\R$. On the other hand,  by the Bochner theorem  for characteristic function,  (1.5) is  characteristic  if and only if $\alpha\in [0, 1]$.

 Here we will extend  Theorem 1 for  complex-valued characteristic functions.  To this end, let us introduce some notions and basic facts.

  If $u$ is a   complex-valued harmonic function  in a domain $D\subset\Co$, then   another complex-valued harmonic function $v$ in $D$ is harmonic conjugate of $u$ if $u+iv$ is analytic in $D$. Recall that if $D\neq\emptyset$, then the harmonic conjugate $v$ of $u$ is unique, up to adding a constant.

Let $f: \R\to\Co$ be a continuous bounded function. The integral
\[
v_f(x,y)= \frac1{\pi}\int_{-\infty}^{\infty} \frac{x-t}{(x-t)^2+y^2} f(t)\,dt
\]
can be chosen as the harmonic conjugate of $u_f$ in $\Cop$  (see, for example, [3, p. 108]). In that case $u_f+iv_f$ coincides with the  usual  Cauchy  transform (the Cauchy integral) of $f$
\begin{equation}\label{1.6}
k_{f}(z)= \frac{i}{\pi}\int_{-\infty}^{\infty} \frac{f(t)}{z-t}\,dt.
\end{equation}
Here the integral is absolutely convergent as long as
\begin{equation}\label{1.7}
\int_{-\infty}^{\infty}\frac{|f(t)|}{1+|t|}\,dt<\infty.
\end{equation}
In particular, this condition is satisfied for any $f\in L^p(\R)$, $1\le p<\infty$, but not in the case of an arbitrary $f\in L^{\infty}(\R)$.  In general, if  $f\in L^p(\R)$, $1\le p\le\infty$, then we shall use the followings harmonic conjugate of (1.2)
\[
V_f(x,y)= \frac1{\pi}\int_{-\infty}^{\infty} \Bigl(\frac{x-t}{(x-t)^2+y^2}+\frac{t}{t^2+1}\Bigr) f(t)\,dt
\]
 (see [3, p.p. 108-109]). Denote
\begin{equation}\label{1.8}
K_{f}(z)= \frac{i}{\pi}\int_{-\infty}^{\infty} \Bigl(\frac{1}{z-t} +\frac{t}{t^2+1}\Bigr) f(t)\,dt,
\end{equation}
$z\in \Co\setminus\R$. We call $K_f$ the  modified   Cauchy  transform of $f$. Both integrals in  (1.6) and in (1.8)   define  analytic functions in  $\Co\setminus\R= \Co^{+}\cup\Co^{-}$ (e.g., see [6, p.p.  144-145]), where $\Co^{-}$ denotes the open lower half-plane in $\Co$. Let us write $k_{f}^{(\pm)}(z)=k_{f}(z)$ for $z\in \Co^{\pm}$, and  $K_{f}^{(\pm)}(z)=K_{f}(z)$ for $z\in \Co^{\pm}$.

  A function $\omega(x)$  is said to be absolutely monotonic  on $(a, b)$ if and only if $\omega(-x)$ is completely monotonic on $(-b, -a)$ (see [7, p.p. 144-145]). It is obvious that such a function $\omega(x)$ can be  characterized by  the inequalities
\begin{equation}\label{1.9}
\omega^{(n)}(x)\ge 0,
\end{equation}
  $n=0,1,2,\dots$.

 The main results are the following theorems:

\begin{thm}.\label{t2}
Suppose that $f: \R\to\Co $ is a bounded continuous  function   and  $f(0)=1$. Then $f$ is a characteristic function if and only if there is a constant $a_f\in\R$ such that: \

({\bf i}) \   $ a_f+K^{(+)}_f(iy)$  is completely monotonic for $y\in (0,\infty)$,
and

({\bf ii}) \   $ -\Bigl(a_f+K^{(-)}_f(iy)\Bigr)$  is absolutely monotonic for $y\in (-\infty, 0)$.
\end{thm}

\begin{thm}.\label{t3}
Suppose that $f$ is as in Theorem 2 and satisfies (1.7). Then $f$ is a characteristic function if and only if:

({\bf i}) \   $ k^{(+)}_f(iy)$  is completely monotonic for $y\in (0,\infty)$,
and

({\bf ii})  \  $ -k^{(-)}_f(iy)$  is absolutely monotonic for $y\in (-\infty, 0)$.
\end{thm}

\section{Proofs}
 \label{s2}

 A  function  $\varphi$ on $\R$ is said to be positive definite if
\[
\sum_{i,j=1}^{n} \varphi(x_i-x_j) c_i\overline{c_j}\ge 0
\]
for every choice of $x_1,\dots,x_n\in\R$, for every choice of $c_1,\dots, c_n\in\Co$, and all $n\in\mathbb{N}$. By the Bochner theorem  (e.g., see [2, Theorem 33.3]), we have that  a   continuous function $\varphi$ is positive definite if and only if there exists  a  non-negative finite measure $\mu$ on $\R$ such that $\varphi(x)=\widehat{\mu}(-x)$.

  We will need to use later  the following lemma:
 \begin{lem}\label{s1}{\rm{[2,  Theorem 33.10].}}
Let  $\varphi$ be a continuous positive definite function on $\R$ such that $\varphi\in L^1(\R)$. Then $\widehat{\varphi}$ is nonnegative, $\widehat{\varphi}$ is in $ L^1(\R)$, and $(\check{\varphi})^{\widehat{}}(x)=\varphi(x)$ for all $x\in\R$.
\end{lem}

 Here
 \[
  \check{\varphi}(t)=\frac1{2\pi}\int_{-\infty}^{\infty}e^{it\xi}\varphi(\xi)\,d\xi
  \]
is  the inverse Fourier transform of  $\varphi$.   We first prove our Theorem 3.

{\bf{ Proof of Theorem 3.}} \ \ 
Suppose that  $f$ is a   characteristic function.  Let $z=x+iy\in \Co$ with $y\neq 0$.   According to (1.7),  the integral in (1.6) converges absolutely. Moreover, it is easily checked [6, p.p. 144-145] that
\begin{equation}\label{2.1}
\frac{d^n}{dy^n}\,k_{f}(iy)= \frac{i}{\pi}\int_{-\infty}^{\infty} \frac{\partial^n}{\partial\,y^n}\Bigl(\frac{1}{iy-t}\Bigr) f(t)\,dt
\end{equation}
for all $m=0,1,2,\dots$.

Let now $z\in\Cop$.  Then by direct calculation we obtain that
\[
\int_0^{\infty} x^n e^{-yx}e^{-ixt}\,dx=(-1)^n\frac{\partial^n}{\partial\,y^n}\Bigl(\frac{i}{iy-t}\Bigr).
\]
This, together with  Bochner's theorem shows that for any $y>0$ the right side of the previous equality is a continuous  positive definite as a function of $t\in\R$.
Set
\[
\varphi_n(t)= (-1)^n\frac{\partial^n}{\partial\,y^n}\Bigl(\frac{i}{iy-t}\Bigr) f(t).
\]
Under the condition (1.7), we have  that $\varphi_n$ satisfies the hypotheses of Lemma 1. Hence
\[
(-1)^n\int_{-\infty}^{\infty} \frac{\partial^n}{\partial\,y^n}\Bigl(\frac{i}{iy-t}\Bigr)f(t)\,dt=\int_{-\infty}^{\infty}\varphi_n(t)\,dt=\widehat{\varphi}_n(0)\ge 0
\]
for all $m=0,1,2,\dots$. Now, by   (2.1),   we get  that the function  $y\to k_{f}^{(+)}(iy)$ satisfies (1.4), i.e., it is completely monotonic on $(0,\infty)$.

In the case where $z\in\Cod$, we have
\[
 \int_{-\infty}^0 |x|^n e^{-yx}e^{-ixt}\,dx=-\frac{\partial^{n}}{\partial\,y^n}\Bigl(\frac{i}{iy-t}\Bigr)
\]
for  $n=0,1,2,\dots$. Again,  we see that the right side of this equality  is continuous and  positive definite as a function of $t\in\R$.   Applying now Lemma 1 to
\[
\varphi_n(t)= -\frac{\partial^n}{\partial\,y^n}\Bigl(\frac{i}{iy-t}\Bigr) f(t),
\]
  we get  as in the previous case that the function   $y\to -k_{f}^{(-)}(iy)$ satisfies  (1.9). Thus,  it is absolutely monotonic on $(-\infty, 0)$.

 Suppose that  $ k_f^{(+)}(iy)$ is completely monotonic for $y\in (0,\infty)$. By the Bernstein-Widder theorem (see [7,  p. 116]), there exists a nonnegative (not necessarily finite) measure $\mu$  supported on $[0,\infty)$ such that
\begin{equation}\label{2.2}
k_f^{(+)}(iy)= \int_0^{\infty} e^{-y t}\,d\mu (t),
\end{equation}
where the integral converges for all $y\in(0,\infty)$. This means that for any  given $\tau>0$, the function
 \begin{equation}\label{2.3}
u_{\tau}(iy):= k_f^{(+)}(i(y+\tau))
\end{equation}
 is completely monotonic for $y\in [0,\infty)$. Now, there is a finite nonnegative measure $\mu_{\tau}$ on $[0,\infty)$ such that
\[
u_{\tau}(iy)=\int_0^{\infty}e^{-y t}\, d\mu_{\tau}(t),
\]
 $y\in [0,\infty)$. Any such $u_{\tau}$ can be  continued analytically to $\Cop$ as the Laplace transform of finite $\mu_{\tau}$.  According to (2.3), we have that the Laplace transform of $\mu$
\begin{equation}\label{2.4}
L_{\mu}(z)=\int_0^{\infty} e^{izt}\,d\mu(t)
\end{equation}
is  well-defined and  is also analytic in  $\Cop$. With equation (2.2) in mind, the applications of  the uniqueness theorem for analytic functions (1.6) and (2.4) in  $\Cop$ yields
\begin{equation}\label{2.5}
\frac{i}{\pi}\int_{-\infty}^{\infty} \frac{f(t)}{z-t}\,dt = k_f(z)= k_f^{(+)}(z)=\int_0^{\infty}e^{izt}\,d\mu(t),
\end{equation}
 $z\in\Cop$.

 Let $y\to -k^{(-)}_f(iy)$  be absolutely monotonic on $ (-\infty, 0)$. By definition, the function  $ -k^{(-)}_f(-iy)$ is completely monotonic for $y\in (0, \infty)$. Using the same argument as before, we have that there exists a nonnegative measure $\eta$  on $(-\infty,0]$ such that
\begin{equation}\label{2.6}
 k_f(z)= k_f^{(-)}(z)=-\int_{-\infty}^0 e^{izt}\,d\eta(t),
\end{equation}
where the integral is  absolutely convergent for each  $z\in \Co^{-}$.

  Fix $y>0$. Using (2.5) and (2.6), we get
\begin{equation}\label{2.7}
\Bigl(P_{y}\ast  f\Bigr)(x)= \frac12\Bigl[k_f(x+iy)-k_f(x-iy)\Bigr]=\int_{-\infty}^{\infty}e^{ixt}\,d\vartheta_y(t),
\end{equation}
where
\[
 \vartheta_y(t)=\frac12e^{-|y|t}\Bigl[\mu(t)+\eta(t)\Bigr].
 \]
Since the integrals (2.4) and (2.6) are absolutely convergent, it follows that $\vartheta_y$ is a finite measure on $\R$.  Now, applying the Bochner theorem to (2.7), we see that $(P_{y}\ast  f)(x)$
is continuous positive definite function on $\R$.  Recall that the family of  Poisson kernels $(P_y)_{y>0}$ form  an approximate unit in $L^1(\R)$ (e.g., see [3,  p. 111]). Hence
\[
\lim_{y\to 0} \Bigl(P_{y}\ast  f\Bigr)(x_0)=f(x_0)
\]
for any point $x=x_0$ of continuity of $f$.  According to the  fact that the pointwise limit of positive definite functions also is a positive definite, we have that $f$ is positive definite on $\R$. This  proves Theorem 3.

{\bf{Proof of Theorem 2}}. \ \
 Fix $z\in\Co\setminus\R$. Since the modified Cauchy kernel
\[
\gamma(z;t)=\frac{i}{\pi}\Bigl(\frac1{z-t}+\frac{t}{t^2+1}\Bigr)
\]
is  an integrable  function, it follows that its Fourier transform is  well-defined and  is also  continuous function. Let $ \zeta_A$ denote  the indicator function of a measurable  subset $A\subset\R$. If $z=iy$,  $y\neq 0$, then by direct calculation we obtain that
\begin{equation}\label{2.8}
\widehat{\gamma}(iy;-x)=\int_{-\infty}^{\infty}\gamma(iy;t) e^{ixt}\,d t= \left\{
\begin{array}{rl}
2\zeta_{[0,\infty)}(x)e^{-yx}-{\rm{sign}}(x)e^{-|x|}, & \mbox{if } x\neq 0,\\
1,  & \mbox{if } x= 0,
\end{array}\right.
\end{equation}
for $y>0$, and
\begin{equation}\label{2.9}
\widehat{\gamma}(iy;-x)=\left\{
\begin{array}{rl}
-2\zeta_{(-\infty,0]}(x) e^{-yx}-{\rm{sign}}(x)e^{-|x|}, & \mbox{if } x\neq 0,\\
-1,  & \mbox{if } x= 0,
\end{array}\right.
\end{equation}
for $y<0$.

 To prove the necessity implication, suppose that $f$ is a   characteristic function. Then
\begin{equation}\label{2.6}
f(t)=\int_{-\infty}^{\infty} e^{ixt}\,d\sigma(x),
\end{equation}
where $\sigma$ is a probability measure on $\R$. Since the function (2.8) is positive for small $y$ and is negative for large $y$,  we have that  $\gamma(iy;t)$   is not necessary  positive definite as a function of $x\in\R$ in general. Therefore, we cannot apply   Lemma 1, as in the proof of Theorem 3. But, on the other hand,  both  (2.8) and (2.9) are   integrable functions for $x\in\R$.

 Let  $y>0$. Substituting (2.10) into (1.8),  applying (2.8) and  using  Fubini's theorem, we get
\[
K_f(iy)=\int_{-\infty}^{\infty}\gamma(iy; t)f(t)\,dt=\int_{-\infty}^{\infty}\widehat{\gamma}(iy; -x)\,d\sigma(x)= \int_0^{\infty} e^{-yx}\,d\sigma_1(x)-a_f,
\]
where $\sigma_1$ denotes the  nonnegative finite measure $2\sigma\cdot \zeta_{[0,\infty)}- \sigma\{0\}$, while
\begin{equation}\label{2.11}
a_f= \int_{-\infty}^{\infty}{\rm{sign}}(x)e^{-|x|}\,d\sigma(x).
\end{equation}
Here $\sigma\{0\}$ is the measure $\sigma$ of the one-point set $\{0\}$. So by the Bernstein-Widder theorem [7,  p. 116], the function  $K_f(iy)+a_f$ is completely monotonic for $y>0$.

 If  $y<0$, then  combining (1.8), (2.9), and (2.10) we obtain
\[
-K_f(iy)= \int_{-\infty}^0 e^{-yx}\,d\sigma_2(x)+a_f
\]
with $\sigma_2= 2\sigma\cdot \zeta_{(-\infty,0]}- \sigma\{0\}$ and   $a_f$ defined in  (2.11).  Finally, it is easy to verify by straightforward calculation of derivatives that the function
\[
-\Bigl(a_f+K^{(-)}_f(iy)\Bigr)=\int_{-\infty}^0 e^{-yx}\,d\sigma_2(x)
\]
satisfies  (1.9) for $y\in (-\infty, 0)$.

  The sufficiency can be proved in a manner similar to the proof of the sufficiency of Theorem 3.

\bibliographystyle{elsarticle-harv}

\begin{thebibliography}{00}

\bibitem[1]{Eg}A.V. Egorov,  On the theory of characteristic functions,  Russian Math. Surveys 59(3) (2004)  567-568.
\bibitem[2]{HR}E. Hewitt,  K.A. Ross,  Abstract Harmonic Analysis, vol. 2,  Springer, Berlin-Heidelberg, 1997.
\bibitem[3]{Ko}P. Koosis,   Introduction to $H_p$  Spaces,  2nd ed.,  Cambridge Tracts in Mathematics,  vol. 115,   Cambridge University Press,  1998.
\bibitem[4]{Lu}E. Lukacs,  Characteristic Functions,  2nd ed.,   Hafner Publishing Co., New York, 1970.
\bibitem[5]{No}S. Norvidas,  On harmonic continuation of characteristic functions,  Lithuanian Math. J.  50(4) (2010)  418--425.
\bibitem[6]{Vl}V.S. Vladimirov,   Methods of the Theory of Generalized Functions,  Taylor$\&$Francis, London,  2002.
\bibitem[7]{Wi}D.V. Widder,  The Laplace Transform,  Princeton Mathematical Series, vol. 6,  Princeton University Press, 1941.
 \end{thebibliography}

{\large{

}}
\end{document}